\documentclass[10pt]{amsart}

\title[Hodge integrals in FJRW theory]{Hodge integrals in FJRW theory}

\author[Gu\'er\'e]{J\'er\'emy Gu\'er\'e}
\address{Humboldt Universit\"at zu Berlin\\
Berlin\\
Deutschland}
\email{jeremy.guere@hu-berlin.de}

\usepackage{amsrefs}
\usepackage{amsfonts}
\usepackage{amsmath}
\usepackage{amssymb}
\usepackage{amscd}
\usepackage{array}
\usepackage{subfigure}
\usepackage{tikz}

\usepackage[all]{xy}

\allowdisplaybreaks[3]

\newcommand{\grj}{\mathfrak{j}}

\newcommand{\kK}{\mathbf{K}}
\newcommand{\fc}{\mathfrak{c}}
\newcommand{\PV}{\mathbf{PV}}

\newcommand{\aA}{\mathbb{A}}
\newcommand{\bB}{\mathfrak{B}}
\newcommand{\rR}{\mathfrak{C}}
\newcommand{\ci}{\mathrm{i}}
\newcommand{\PP}{\mathbb{P}}
\newcommand{\CC}{\mathbb{C}}

\newcommand{\cO}{\mathcal{O}}

\renewcommand{\(}{\left(}
\renewcommand{\)}{\right)}

\newcommand{\ii}{\mathrm{i}}
\newcommand{\cC}{\mathcal{C}}

\newcommand{\cL}{\mathcal{L}}
\newcommand{\cM}{\mathcal{M}}

\newcommand{\cQ}{\mathcal{Q}}

\newcommand{\fq}{\mathfrak{q}}

\newcommand{\cS}{\mathrm{Sym}}
\newcommand{\sS}{\mathcal{S}}

\newcommand{\Ch}{\mathrm{Ch}}
\newcommand{\Td}{\mathrm{Td}}

\newcommand{\st}{\mathbf{H}}

\newcommand{\cvir}{c_{\mathrm{vir}}}
\newcommand{\cvirPV}{c_{\mathrm{vir}}^{\mathrm{PV}}}

\DeclareMathOperator{\Aut}{Aut}

\theoremstyle{plain}

\newtheorem{thm}{Theorem}[section]

\newtheorem*{lem*}{Lemma}

\newtheorem{cor}[thm]{Corollary}

\theoremstyle{definition}

\newtheorem{dfn}[thm]{Definition}

\newtheorem{rem}[thm]{Remark}
\newtheorem*{rem*}{Remark}
\newtheorem*{rems*}{Remarks}
\newtheorem{exa}[thm]{Example}
\newtheorem*{exa*}{Example}
\newtheorem*{exait*}{\rm \em Example}
\newtheorem*{exadefit*}{\rm \em Example/Definition}
\newtheorem*{cla*}{\rm \em Claim}

\newtheorem*{dfn*}{Definition}

\usepackage{amsxtra}

\def\<{\left\langle}
\def\>{\right\rangle}

\begin{document}
\begin{abstract}
We study higher genus Fan--Jarvis--Ruan--Witten theory of any chain polynomial with any group of symmetries.
Precisely, we give an explicit way to compute the cup product of Polishchuk and Vaintrob's virtual class with the top Chern class of the Hodge bundle.
Our formula for this product holds in any genus and without any assumption on the semi-simplicity of the underlying cohomological field theory.
\end{abstract}

\maketitle

\tableofcontents

\setcounter{section}{-1}

\section{Introduction}
In $1999$, Candelas, de la Ossa, Green, and Parkes \cite{Candelas} proposed a famous formula for the genus-zero invariants enumerating rational curves on the quintic threefold.
It has later been proved by Givental \cite{Giv0,Gi} and Lian--Liu--Yau \cite{LLY,LLY2,LLY3}, giving a full understanding of Gromov--Witten invariants in genus zero for the quintic threefold.
The genus-one case was then completely solved by Zinger \cite{Zi}. 
However, we still lack a complete understanding in higher genus.

In fact, even the problem of computing genus-zero Gromov--Witten invariants of projective varieties is not completely solved. One of the techniques is called quantum Lefschetz principle (see for instance \cite{Coates3}) and compares Gromov--Witten invariants of a complete intersection with those of the ambient projective space. Thus, we are still missing Gromov--Witten invariants attached to primitive cohomological classes, i.e.~the classes which do not come from the ambient space.

When considering complete intersections in weighted projective spaces, the theory for genus-zero and with ambient cohomological classes looks as complicated as in the higher genus case, because of the lack of a convenient assumption: convexity.
Convexity hypothesis roughly turns the virtual fundamental cycle from Gromov--Witten theory into the top Chern class of a vector bundle, making it easier to compute. But in general, this assumption is not satisfied and the quantum Lefschetz principle can fail \cite{Coates1}.

In \cites{FJRW,FJRW2}, Fan, Jarvis, and Ruan, based on ideas of Witten \cite{Witten}, have switched to another quantum theory which they define for polynomial singularities. We call it FJRW theory and it is attached to a Landau--Ginzburg orbifold $(W,G)$, where $W$ is a non-degenerate quasi-homogeneous polynomial singularity and $G$ is a group of diagonal symmetries of $W$.
The Landau--Ginzburg/Calabi--Yau correspondence conjecture \cite{Chiodo2} describes, under some Calabi--Yau assumption, the relation between this new theory and Gromov--Witten theory of the hypersurface\footnote{More precisely, the Landau--Ginzburg orbifold $(W,G)$ corresponds to the quotient stack $\left[X/\tilde{G}\right]$ where $X$ is the hypersurface corresponding to the zero locus of $W$, $\tilde{G}$ is the group $G/\langle\grj\rangle$ and $\grj$ is the matrix defined in \eqref{gradingelement}.} defined by $W$ in the corresponding weighted projective space.
In genus zero, this conjecture has been proven in some convex cases in \cite{LG/CY}.

Therefore, the study of FJRW theory appears as a new point of view toward the study of Gromov--Witten theory.
In \cite{Guere1}, we described an explicit way to compute FJRW theory in genus zero for polynomials whose Gromov--Witten counterparts are unknown, because of the lack of convexity\footnote{The corresponding notion in FJRW theory is called concavity.}.
In the recent work \cite{Li15}, the Landau--Ginzburg/Calabi--Yau correspondence is studied in higher genus for the quintic hypersurface in $\PP^4$.

In non-zero genus, both Gromov--Witten and Fan--Jarvis--Ruan--Witten theories are extremely difficult to compute.
There are nevertheless some powerful techniques, as the localization \cite{local} and the degeneration \cite{degen} formulas in Gromov--Witten theory and Teleman's reconstruction theorem for conformal generically semi-simple cohomological field theories \cite{Teleman}.
For instance, the localization formula determines all Gromov--Witten invariants of homogeneous spaces \cites{Kontdef,local}.
Also, Teleman's reconstruction theorem takes a major place in the proof of the generalization of Witten conjecture to ADE singularities \cites{FSZ,FJRW}, in the proof of Pixton's relations \cite{Pixton}, and more recently in the study of higher-genus mirror symmetry \cite{impost} after Costello--Li \cite{Costel}.

The method presented in this paper is quite different from above techniques and is valid for a range of Landau--Ginzburg orbifolds for which no previous techniques are applicable.
More precisely, it works without any semi-simplicity assumption and it uses instead the K-theoretic vanishing properties of a recursive complex of vector bundles. It is a direct generalization of the results in \cite{Guere1}, where recursive complexes are introduced for the first time.

In this introduction, we state our theorem in the chain case with the so-called narrow condition and refer to Theorem \ref{chainhigh} for a complete statement.
Let $(W,G)$ be a Landau--Ginzburg orbifold, where $W$ is a chain polynomial
\begin{equation*}
W = x_1^{a_1}x_2+\dotsb+x_{N-1}^{a_{N-1}}x_N+x_N^{a_N}
\end{equation*}
and $G$ is a group of diagonal matrices preserving $W$ and containing the matrix $\grj$ defined in \eqref{gradingelement}.

We take $n$ diagonal matrices $\gamma(1),\dotsc,\gamma(n)$ in the group $G$ with no entries equal to $1$ (narrow condition) and we consider the moduli space $\sS_{g,n}(W,G)(\gamma(1),\dotsc,\gamma(n))$ of genus-$g$ $(W,G)$-spin marked curves with monodromy $\gamma(i)$ at the $i$-th marked point (see Section \ref{section1.1} for definitions).

We denote by $\cL_1,\dotsc,\cL_N$ the universal line bundles associated to the variables $x_1,\dotsc,x_N$, by $\cvirPV(\gamma(1),\dotsc,\gamma(n))_{g,n}$ the associated virtual class defined by Polishchuk and Vaintrob \cite{Polish1}, by $\pi$ the morphism from the universal curve to the moduli space, and by $\mathbb{E}:=\pi_*\omega$ the Hodge vector bundle on the moduli space 
corresponding to global differential forms on the curves.

\begin{thm}\label{chainhigh0}
Let $(W,G)$, $\gamma(1),\dotsc,\gamma(n)$, and $\cL_1,\dotsc,\cL_N$ be as above.
For any genus $g$, we have
\begin{equation*}
c_\mathrm{top}(\mathbb{E}^\vee) ~ \cvirPV(\gamma(1),\dotsc,\gamma(n))_{g,n}  =  \lim_{t_1 \to 1}  \prod_{j=1}^N \fc_{t_j}(-R^\bullet \pi_*(\cL_j)) \cdot \fc_{t_{N+1}}(\mathbb{E}^\vee),
\end{equation*}
where the variables $t_j$ satisfy the relations $t_j^{a_j} t_{j+1} = 1$ for $j \leq N$, and where the function $\fc_t$ is the characteristic class introduced in \cite{Guere1}, see also equation \eqref{newclass}.
\end{thm}

Above theorem has several consequences:
\begin{itemize}
\item computation of Hodge integrals in FJRW theory via a computer \cites{PhDJG,computerprogram},
\item computation of double ramification hierarchies \cites{Bur14,Guere3},
\item new method to study non-semi-simple cohomological field theories,
\item tautological relations in the Chow ring of the moduli spaces of $(W,G)$-spin curves, in particular in the moduli spaces of $r$-spin and of stable curves \cite{team4}.
\end{itemize}

\begin{rem}
It is important to strike the fact that FJRW theory is not a generically semi-simple cohomological field theory in general, especially when the group $G$ is generated by the element $\grj$ defined in \eqref{gradingelement}. In such cases, Teleman's reconstruction theorem \cite{Teleman} does not apply and, to our knowledge, the method presented in this paper is the first comprehensive approach in higher genus for these theories, although we only obtain a partial information on the virtual class.
\end{rem}

Integrals of the form
\begin{equation*}
\int_{\overline{\cM}_{g,n}} c_\mathrm{top}(\mathbb{E}^\vee) ~ \alpha ~ , ~~ \alpha \in H^*(\overline{\cM}_{g,n})
\end{equation*}
are called Hodge integrals.
Thus Theorem \ref{chainhigh0}, together with Mumford's \cite{Mumford} and Chiodo's \cite{Chiodo1} formulas, yields an explicit way to compute Hodge integrals in FJRW theory in any genus, and it has been implemented into a computer, see \cites{PhDJG,computerprogram}.
In particular, it is used in \cite{Guere3} to provide a positive answer to Buryak's conjecture \cite{Bur14} on the double ramification hierarchy for $r$-spin theory with $r \leq 5$.
\begin{thm}[See {\cite[Theorem 1.1]{Guere3}}]
For the $3$-spin theory the double ramification hierarchy coincides with the Dubrovin--Zhang hierarchy. For the $4$ and $5$-spin theories the double ramification hierarchy is related to the Dubrovin--Zhang hierarchy by the following Miura transformation:
\begin{align*}
&\left\{
\begin{aligned}
&w^1=u^1+\frac{\epsilon^2}{96}u^3_{xx},\\
&w^2=u^2,\\
&w^3=u^3,
\end{aligned}\right.&&\text{for $r=4$};\\
&\left\{
\begin{aligned}
&w^1=u^1+\frac{\epsilon^2}{60}u^3_{xx},\\
&w^2=u^2+\frac{\epsilon^2}{60}u^4_{xx},\\
&w^3=u^3,\\
&w^4=u^4,
\end{aligned}\right.&&\text{for $r=5$}.
\end{align*}
\end{thm}

Theorem \ref{chainhigh0} has another remarkable consequence: it provides tautological relations in the Chow ring of the moduli space of $(W,G)$-spin curves.
Indeed, the result of Theorem \ref{chainhigh0} holds in the Chow ring and not only in the cohomology ring. Furthermore, it is a statement on the moduli space of $(W,G)$-spin curves, obtained before forgetting the spin structure to end in the moduli space of stable curves.
Even in the $r$-spin case, where the underlying cohomological field theory is generically semi-simple and conformal, these results are new. The main reason is that Teleman's reconstruction theorem \cite{Teleman} only holds in the cohomology ring and after pushing-forward to the moduli space of stable curves.

\begin{cor}
Let $(W,G)$, $\gamma(1),\dotsc,\gamma(n)$, and $\cL_1,\dotsc,\cL_N$ be as above.
For any genus $g$, the expression
\begin{equation*}
\prod_{j=1}^N \fc_{t_j}(-R^\bullet \pi_*(\cL_j)) \cdot \fc_{t_{N+1}}(\mathbb{E}^\vee)
\end{equation*}
from Theorem \ref{chainhigh0} is a Laurent power series in the variable $\epsilon:=t_1^{-1}-1$ of the form
\begin{equation*}
C_{-p} \cdot \frac{1}{\epsilon^p} + C_{-p+1} \cdot \frac{1}{\epsilon^{p-1}} + \dotsb + C_{-1} \cdot \frac{1}{\epsilon} + C_0 + C_1 \cdot \epsilon + \dotsb ,
\end{equation*}
where
\begin{equation*}
C_m \in \bigoplus_{k \geq \mathrm{degvir}+g -m} A^{k}(\sS^G_{g,n})
\end{equation*}
and $p=2g-3+n-\mathrm{degvir}$, the integer $\mathrm{degvir}$ being the Chow degree of $\cvirPV(\gamma(1),\dotsc,\gamma(n))_{g,n}$.
Thus, we obtain tautological relations
\begin{equation*}
C_m = 0 \textrm{ for all } m<0
\end{equation*}
in the Chow ring of the moduli space of $(W,G)$-spin curves.
\end{cor}

It is a work in progress \cite{team4} to compare the push-forward of these relations to the moduli space of stable curves with other tautological relations, e.g.~Pixton's relations \cite{Pixton}.

\noindent
\textbf{Structure of the paper}
In the first part, we briefly recall the main definitions and constructions in FJRW theory following our previous article \cite{Guere1}.
The second part consists of the main theorem \ref{chainhigh} together with its proof.


\noindent
\textbf{Acknowledgement.}
The author is grateful to Alessandro Chiodo and Yongbin Ruan for motivating discussions on this subject.
He also thanks Alexander Buryak, Felix Janda, Paolo Rossi, Dimitri Zvonkine, and Rahul Pandharipande for discussions related to this paper.
He is extremely grateful to Alexander Polishchuk for his comments on the main theorem and its proof.
The author was supported by the Einstein Stiftung.


\section{Quantum singularity theory}\label{section1}
In this section, we give a brief summary of the necessary definitions for the quantum singularity (or FJRW) theory of a Landau--Ginzburg (LG) orbifold.
We use notations of \cite{Guere1} where we dealt with invertible polynomials, but here we are mainly interested in chain or loop polynomials.

\subsection{Conventions and notations}\label{QST}\label{term}
The quantum singularity theory was first introduced by Fan--Jarvis--Ruan \cite{FJRW,FJRW2} after ideas of Witten \cite{Witten}. In particular, Fan, Jarvis, and Ruan constructed a cohomological class called virtual class, via an analytic construction from Witten's initial sketched idea \cite{Witten} formalized for A-singularities by Mochizuki.
Polishchuk and Vaintrob \cite{Polish1} provided an algebraic construction, which generalized their previous construction and that of Chiodo \cite{ChiodoJAG} in the A-singularity case.

We do not know in general whether the two constructions coincide. In FJRW terminology, there is a decomposition of the state space into narrow and broad states.
Chang, Li, and Li prove in \cite[Theorem 1.2]{Li2} the match when only narrow entries occur.
For almost all LG orbifolds $(W,G)$ where $W$ is an invertible polynomial and $G$ is the maximal group of symmetries, we proved in \cite[Theorem 3.25]{Guere1} that the two classes are the same up to a reparametrization of the broad states.
Nevertheless, for smaller groups or more general polynomials, we still do not know whether these two classes coincide.
Therefore, in the whole paper, by virtual class we mean Polishchuk--Vaintrob's version, as soon as we are working with broad states together with non-maximal group $G$.

Furthermore, we work in the algebraic category and over $\CC$.
All stacks are proper Deligne--Mumford stacks; we use also the term ``orbifold'' for this type of stacks.
We denote orbifolds by curly letters, e.g.~$\cC$ is an orbifold curve and the scheme $C$ is its coarse space.
We recall that vector bundles are coherent locally free sheaves and that the symmetric power of a two-term complex is the complex
\begin{equation*}
\cS^k\([ A \rightarrow B ]\)  =  [ \cS^k A \rightarrow \cS^{k-1} A \otimes B \rightarrow \dotso \rightarrow A \otimes \Lambda^{k-1} B \rightarrow \Lambda^k B ]
\end{equation*}
with morphisms induced by $A \rightarrow B$.

All along the text, the index $i$ varies from $1$ to $n$ and refers exclusively to the marked points of a curve whereas the index $j$ varies from $1$ to $N$ and corresponds to the variables of the polynomial.
We represent tuples by overlined notations, e.g.~$\overline{\gamma}=(\gamma(1),\dotsc,\gamma(n))$, or by underlined notations, e.g.~$\underline{p}=(p_1,\dotsc,p_N)$.

\subsection{Landau--Ginzburg orbifold}\label{section1.1}
Let $w_1,\dotsc,w_N$ be coprime positive integers, $d$ be a positive integer and $\fq_j := w_j/d$ for all $j$.
We consider a quasi-homogeneous polynomial $W$ of degree $d$ with weights\footnote{We assume that a choice of coprime positive weights $w_1,\dotsc,w_N$ is unique.} $w_1,\dotsc,w_N$, and with an isolated singularity at the origin. We say that such a polynomial $W$ is non-degenerate.
In particular, for any $\lambda,x_1,\dotsc,x_N \in \CC$, we have
\begin{equation*}
W(\lambda^{w_1} x_1,\dotsc,\lambda^{w_N} x_N) = \lambda^d  W(x_1,\dotsc,x_N)
\end{equation*}
and the dimension of the Jacobian ring
\begin{equation*}
\cQ_W := \CC \left[ x_1,\dotsc,x_N \right] / \left( \partial_1 W,\dotsc, \partial_N W \right)
\end{equation*}
is finite over $\CC$.

An admissible group of symmetries for the polynomial $W$ is a group $G$ made of diagonal matrices $\textrm{diag}(\lambda_1,\dotsc,\lambda_N)$ satisfying
\begin{equation*}
W(\lambda_1 x_1,\dotsc,\lambda_N x_N)=W(x_1,\dotsc,x_N) \quad \textrm{for every } (x_1,\dotsc,x_N) \in \CC^N
\end{equation*}
and containing the grading element
\begin{equation}\label{gradingelement}
\grj := \textrm{diag}(e^{2 \ii \pi \fq_1},\dotsc,e^{2 \ii \pi \fq_N})~, \quad \fq_j := \frac{w_j}{d}.
\end{equation}
The group $G$ is finite and it contains the cyclic group $\mu_d$ of order $d$ generated by $\grj$.
We denote the biggest admissible group by $\textrm{Aut}(W)$.

\begin{dfn}
A Landau--Ginzburg (LG) orbifold is a pair $(W,G)$ with $W$ a non-degenerate (quasi-homogeneous) polynomial and $G$ an admissible group.
\end{dfn}

The quantum singularity theory developed by Fan, Jarvis, and Ruan \cite{FJRW,FJRW2} is defined for any LG orbifold.
In fact, it mostly depends on the weights, the degree, and the group. Precisely, by \cite[Theorem 4.1.8.9]{FJRW}, the theories for two LG orbifolds $(W_1,G)$ and $(W_2,G)$ where the polynomials $W_1$ and $W_2$ have the same weights and degree are isomorphic.

In the context of mirror symmetry, a well-behaved class of polynomials has been introduced by Berglund--H\"ubsch \cite{Hubsch}.
We say that a polynomial is invertible when it is non-degenerate with as many variables as monomials.
According to Kreuzer--Skarke \cite{Kreuzer}, every invertible polynomial is a Thom--Sebastiani (TS) sum of invertible polynomials, with disjoint sets of variables, of the following three types
\begin{equation}\label{ThomSebastiani}
\begin{array}{lll}
\textrm{Fermat:} & \qquad x^{a+1} & \\
\textrm{chain of length } c: & \qquad x_1^{a_1}x_2+\dotsb+x_{c-1}^{a_{c-1}} x_c+x_c^{a_c+1} & (c \geq 2), \\
\textrm{loop of length } l: & \qquad x_1^{a_1}x_2+\dotsb+x_{l-1}^{a_{l-1}} x_l+x_l^{a_l}x_1 & (l \geq 2). \\
\end{array}
\end{equation}

\begin{rem}
In this paper, we consider only polynomials which are of the three types above, and not a Thom--Sebastiani sum of them.
\end{rem}



For any $\gamma \in \mathrm{Aut}(W)$, the set of broad variables with respect to $\gamma$ is
\begin{equation}\label{bgamma}
\bB_\gamma = \left\lbrace x_j \left| \right. \gamma_j =1 \right\rbrace.
\end{equation}

\begin{dfn}\label{statespace}
The state space\footnote{We refer to \cite[Equation (4)]{LG/CY} or \cite[Equation (5.12)]{Polish1} for details about the bidegree and the pairing in this space.} for the LG orbifold $(W,G)$ is the vector space
\begin{eqnarray*}
\st_{(W,G)} & = & \bigoplus_{\gamma \in G} \st_\gamma \\
& = & \bigoplus_{\gamma \in G} (\cQ_{W_\gamma} \otimes d\underline{x}_\gamma)^G,
\end{eqnarray*}
where $W_\gamma$ is the $\gamma$-invariant part of the polynomial $W$, $\cQ_{W_\gamma}$ is its Jacobian ring, the differential form $d\underline{x}_\gamma$ is $\bigwedge_{x_j \in \bB_\gamma} dx_j$, and the upper-script $G$ stands for the invariant part under the group $G$.
\end{dfn}

At last, the quantum singularity theory for an LG orbifold $(W,G)$ is a cohomological field theory, i.e.~the data of multilinear maps
\begin{equation*}
c_{g,n} \colon \st^{\otimes n} \rightarrow H^*(\overline{\cM}_{g,n}),
\end{equation*}
which are compatible under gluing and forgetting-one-point morphisms.
%
More precisely, the maps $c_{g,n}$ factor through the cohomology (and even the Chow ring) of another moduli space $\sS_{g,n}(W,G)$ attached to the LG orbifold $(W,G)$; the map
\begin{equation*}
(\cvir)_{g,n} \colon \st^{\otimes n} \rightarrow A^*(\sS_{g,n}(W,G))
\end{equation*}
is called the virtual class\footnote{For the polynomial $x^r$, we obtain the moduli space of $r$-spin structures and the virtual class is called Witten $r$-spin class.}, where $A^*$ can stand for the cohomology or the Chow ring.
Then, via the natural forgetful morphism $\textrm{o} \colon \sS_{g,n}(W,G) \rightarrow \overline{\cM}_{g,n}$, we get
\begin{equation}\label{pushfor}
c_{g,n} := (-1)^\mathrm{degvir}\frac{\textrm{card}(G)^g}{\mathrm{deg} (\textrm{o})} \cdot \textrm{o}_*(\cvir)_{g,n},
\end{equation}
where $(-1)^\mathrm{degvir}$ acts as $(-1)^m$ on $A^m(\overline{\cM}_{g,n})$.

\begin{rem}
In the case of $r$-spin curves, the degree of the forgetful morphism $\textrm{o}$ above equals $r^{2g-1}$.
In general, for the maximal group $G=\textrm{Aut}(W)$, this degree also equals $r^{2g-1}$, where $r$ is the exponent of the group.
\end{rem}

The moduli space $\sS_{g,n}(W,G)$ is defined in \cite[Section 2]{FJRW} as follows.
First, let us fix $r$ to be the exponent of the group $G$, i.e.~ the smallest integer $l$ such that $\gamma^l=1$ for every element $\gamma \in G$. We recall that an $r$-stable curve is a smoothable\footnote{Concretely, smoothable means that the local picture at the node is $\left[ \left\lbrace xy=0 \right\rbrace / \mu_r \right]$ with the balanced action $\zeta_r \cdot (x,y) = (\zeta_r x , \zeta_r^{-1} y)$.} orbifold curve with markings whose non-trivial stabilizers have fixed order $r$ and are only at the nodes and at the markings. Moreover, its coarse space is a stable curve.

Then, the moduli space $\sS_{g,n}(W,G)$ classifies all $r$-stable curves of genus $g$ with $n$ marked points, together with $N$ line bundles and $s$ isomorphisms
\begin{equation*}
(\cC; \sigma_1,\dotsc,\sigma_n;\cL_1,\dotsc,\cL_N;\phi_1,\dotsc,\phi_s),
\end{equation*}
where the isomorphisms $\phi_1,\dotsc,\phi_s$ give some constraints (see below) on the choice of $\cL_1,\dotsc,\cL_N$.
We call such data a $(W,G)$-spin curve.

To get the constraints $\phi_1,\dotsc,\phi_s$, first choose a Laurent polynomial $Z$ with weights $w_1,\dotsc,w_N$ and degree $d$ just as $W$, with the additional property $\Aut(W+Z)=G$ about the maximal group, see \cite[Prop.~3.4]{Krawitz2}.
Then, denoting by $M_1,\dotsc,M_s$ all the monomials of $W+Z$, we have
\begin{equation}\label{constraints}
\phi_k \colon M_k(\cL_1,\dotsc,\cL_N) \simeq \omega_{\textrm{log}} := \omega_\cC (\sigma_1+\dotsb+\sigma_n)~,~~ \textrm{for all $k$.}
\end{equation}
The moduli space that we obtain does not depend on the choice of the Laurent polynomial $Z$, see \cite{LG/CYquintique}.

A line bundle on an orbifold point comes with an action of the isotropy group at that point, i.e.~locally at a marked point $\sigma_i$, we have an action
\begin{equation}\label{multiplicities}
\zeta_r \cdot (x,\xi) = (\zeta_r x, \zeta_r^{m_j(i)} \xi) ~, \quad \textrm{with } m_j(i) \in \left\lbrace 0 , \dotsc, r-1 \right\rbrace
\end{equation}
called the monodromy of the line bundle $\cL_j$ at the marked point $\sigma_i$.
Since the logarithmic canonical line bundle $\omega_{\textrm{log}}$ is a pull-back from the coarse curve, then its multiplicity is trivial on each marked point, so that equations \eqref{constraints} give
\begin{equation*}
\gamma(i):=(e^{2 \ci \pi m_1(i)/r},\dotsc,e^{2 \ci \pi m_N(i)/r}) \in \Aut(W+Z) = G.
\end{equation*}
We define the type of a $(W,G)$-spin curve as $\overline{\gamma}:=(\gamma(1),\dotsc,\gamma(n)) \in G^n$.
It yields a decomposition
\begin{equation*}
\sS_{g,n}(W,G) = \bigsqcup_{\overline{\gamma} \in G^n} \sS_{g,n}(W,G)(\gamma(1),\dotsc,\gamma(n)),
\end{equation*}
where $\sS_{g,n}(W,G)(\overline{\gamma})$ is an empty component when the selection rule
\begin{equation}\label{selecrule}
\gamma(1) \dotsm \gamma(n) = \grj^{2g-2+n}
\end{equation}
is not satisfied, see \cite[Proposition 2.2.8]{FJRW}.

\subsection{$\Aut(W)$-invariant states}
From the definition \ref{statespace} of the state space $\st_{(W,G)}$, we see that it always contain the subspace
\begin{eqnarray*}
\st_{(W,G),\mathrm{Aut}(W)} & = & \bigoplus_{\gamma \in G} (\cQ_{W_\gamma} \otimes d\underline{x}_\gamma)^{\mathrm{Aut}(W)} \\
& \subset & \st_{(W,G)}.
\end{eqnarray*}

\begin{dfn}
The subspace $\st_{(W,G),\mathrm{Aut}(W)}$ is called the $\mathrm{Aut}(W)$-invariant part.
\end{dfn}

For invertible polynomials $W$ and any group $G$, the $\Aut(W)$-invariant part has a particularly nice and explicit description, see \cites{Krawitz,Guere1}.
Using the language from \cite{Guere1}, we can attach a graph $\Gamma_W$ to any invertible polynomial, illustrating its Kreuzer--Skarke decomposition as a Thom--Sebastiani sum of Fermat, chain, and loop polynomials.
Then, we consider decorations\footnote{A decoration $\rR_\gamma$ is a subset of the set of broad variables $\bB_\gamma = \left\lbrace x_j \left| \right. \gamma_j =1 \right\rbrace$. Definitions for admissible and balanced are in \cite[Definition 1.5]{Guere1}.} $\rR_\gamma$ of the graph $\Gamma_W$ that are admissible and balanced, and to each such decoration we associate an explicit element $e(\rR_\gamma)$ of $\st_{(W,G),\mathrm{Aut}(W)}$.
At last, by \cite{Krawitz} and \cite[Equation (10)]{Guere1}, the set of all these elements forms a basis of $\st_{(W,G),\mathrm{Aut}(W)}$.

\begin{exa}
Let $W=x_1^{a_1}x_2+\dotsb+x_{c-1}^{a_{N-1}} x_N+x_N^{a_N+1}$ be a chain polynomial.
For any element $\gamma \in G$, the set of broad variables is of the form
\begin{equation*}
\bB_\gamma = \left\lbrace x_{b+1}, \dotsc,x_N \right\rbrace
\end{equation*}
and there is exactly one admissible decoration $\rR_\gamma$ given by
\begin{equation*}
\rR_\gamma = \left\lbrace x_{N-2j} \left| \right. N-2j > b \right\rbrace.
\end{equation*}
This decoration is balanced if and only if $N-b$ is even and the corresponding element is
\begin{equation*}
e_\gamma := e(\rR_\gamma) = \biggl( \prod_{\substack{b < j \leq N \\ N-j \textrm{ odd}}} a_j x_j^{a_j-1}\biggr)  \cdot d x_{b+1} \wedge \dotsm \wedge d x_N.
\end{equation*}
\end{exa}

\begin{exa}
Let $W=x_1^{a_1}x_2+\dotsb+x_{c-1}^{a_{N-1}} x_N+x_N^{a_N}x_1$ be a loop polynomial.
For an element $1 \neq \gamma \in G$, the set of broad variables is empty.
For the identity element, it is $\bB_1 = \left\lbrace x_1, \dotsc,x_N \right\rbrace$.
Then, if $N$ is odd, there is no admissible and balanced decoration.
But if $N$ is even, we have two distinct admissible and balanced decorations given by
\begin{equation*}
\rR^+_1 = \left\lbrace x_j \left| \right. j \textrm{ even} \right\rbrace \textrm{ and } \rR^-_1 = \left\lbrace x_j \left| \right. j \textrm{ odd} \right\rbrace,
\end{equation*}
and the two corresponding elements are
\begin{equation*}
e^+ := e(\rR^+_1) = \biggl( \prod_{x_j \textrm{ odd}} a_j x_j^{a_j-1} - \prod_{x_j \textrm{ even}} -x_j^{a_j-1}\biggr)  \cdot d x_1 \wedge \dotsm \wedge d x_N,
\end{equation*}
and $e^-$ exchanging even and odd.
\end{exa}

\subsection{Sketch of the definition of PV virtual class}
Polishchuk--Vaintrob construction \cite{Polish1} of the virtual class $(\cvir)_{g,n}$ for an LG orbifold $(W,G)$ uses the notion of matrix factorizations. We briefly recall the main steps.

Consider a component $\sS_{g,n}(\overline{\gamma})$ of type $\overline{\gamma}=(\gamma(1),\dotsc,\gamma(n)) \in G^n$.
We denote by $\pi$ the projection of the universal curve to this component and we look at the higher push-forwards $R^\bullet\pi_*\cL_j$ of the universal line bundles.
We take resolutions of $R^\bullet\pi_*\cL_j$ by complexes $[A_j \rightarrow B_j]$ of vector bundles and we set
\begin{equation*}
X := \mathrm{Spec} ~ \mathrm{Sym} \bigoplus_{j=1}^N A_j^\vee \quad \textrm{and} \quad p \colon X \rightarrow \sS_{g,n}(\overline{\gamma}).
\end{equation*}
The differential $[A_j \rightarrow B_j]$ induces a section $\beta$ of the vector bundle $p^*\bigoplus_j B_j$ on $X$.
Polishchuk and Vaintrob show how to construct a section $\alpha$ of the dual vector bundle $p^*\bigoplus_j B_j^\vee$, using the algebraic relations \eqref{constraints} between the line bundles $\cL_1,\dotsc,\cL_N$.
The choice of the resolutions and the existence of the section $\alpha$ require several steps, see \cite[Section 4.2 Steps 1-4]{Polish1}.

Using evaluation of the line bundles $\cL_1,\dotsc,\cL_N$ at the marked points, they also construct a morphism
\begin{equation}\label{aAffinespace}
Z \colon X \rightarrow \aA^{\overline{\gamma}}:= \prod_{i=1}^n \left( \aA^N \right)^{\gamma(i)},
\end{equation}
where $\left( \cdot \right)^{\gamma(i)}$ is the fixed locus under the action of $\gamma(i)$.
In particular, the set of coordinates of the affine space $\aA^{\overline{\gamma}}$ is indexed as
\begin{equation*}
\left\lbrace x_j(i)\right\rbrace_{(\sigma_i,x_j) \in \bB_{\overline{\gamma}}} ~,~~ \textrm{where  }\bB_{\overline{\gamma}} = \left\lbrace (\sigma_i,x_j) \left| \right. \gamma_j(i) = 1 \right\rbrace
\end{equation*}
and we further consider the invertible polynomial $W_{\overline{\gamma}}$ on $\aA^{\overline{\gamma}}$ given by
\begin{equation*}
W_{\overline{\gamma}} := W_{\gamma(1)}(x_1(1),\dotsc,x_N(1)) + \dotsb + W_{\gamma(n)}(x_1(n),\dotsc,x_N(n)),
\end{equation*}
where $W_{\gamma(i)}$ is the restriction of $W$ to $\left( \aA^N \right)^{\gamma(i)}$.

At last, the two sections $\alpha$ and $\beta$ yield a Koszul matrix factorization $\PV$ on $X$.
Polishchuk and Vaintrob checked that the potential of $\PV$ is precisely the function $-Z^*W_{\overline{\gamma}}$ on $X$.
To sum up, we have
\begin{center}
\begin{tikzpicture}
\node (X) at (0.8,1) {$X$};
\node[above] (E) at (1.8,1.6) {$\quad \quad \quad \quad \quad \quad \quad  \PV \in \textrm{MF}(X,-Z^*W_{\overline{\gamma}})$};
\node (A) at (0,0) {$\aA^{\overline{\gamma}}$};
\node (S) at (1.6,0) {$S$};
\draw[->,>=stealth] (E) to[bend left=10] (X);
\draw[->,] (X) -- (A);
\draw[->] (X) -- (S);
\draw (0.45,0.65) node[left] {$Z$};
\draw (1.15,0.65) node[right] {$p$};
\end{tikzpicture}
\end{center}

The matrix factorization $\PV$ is used as a kernel in a Fourier--Muka\"i transform
\begin{equation}\label{foncteur}
\begin{array}{lcccc}
\Phi \colon & \textrm{MF}(\aA^{\overline{\gamma}},W_{\overline{\gamma}}) & \longrightarrow & \textrm{MF}(S,0) \\
            &                               U                            & \longmapsto     & p_*(Z^*(U) \otimes \PV),   \\
\end{array}
\end{equation}
where the two-periodic complex $Z^*(U) \otimes \PV$ is supported inside the zero section $S \hookrightarrow X$ (see \cite[Sect.~4.2, Step 4; Proposition 1.4.2]{Polish1}), so that the push-forward functor is well-defined.

Polishchuk and Vaintrob proved that the Hochschild homology of the category of matrix factorizations on an affine space with polynomial potential $f(y_1,\dotsc,y_m)$ is isomorphic to $\cQ_f \otimes dy_1 \wedge \dotsc \wedge dy_m$. They also give a very explicit description of the Chern character map.
We then have the commutative diagram
\begin{equation*}
\xymatrix{
    \mathrm{MF}(\aA^{\overline{\gamma}},W_{\overline{\gamma}}) \ar[r]^{\Phi} \ar[d]_\Ch & \mathrm{MF}(S,0) \ar[d]^\Ch \\
    \otimes_{i=1}^n \st_{\gamma(i)} \ar[r]_{\Phi_*} & H^*(S)
    }
\end{equation*}
At last, given states $u_{\gamma(1)}, \dotsc, u_{\gamma(n)}$ such that $u_{\gamma(i)} \in \st_{\gamma(i)}$, the virtual class evaluated at these states is
\begin{equation}\label{virtualcomplex}
(\cvir)_{g,n}(u_{\gamma(1)}, \dotsc, u_{\gamma(n)}) = \Phi_*(u_{\gamma(1)}, \dotsc, u_{\gamma(n)}) ~ \prod_{j=1}^N \frac{\Td (B_j)}{\Td (A_j)}
\end{equation}
and is an element of $A^*(\sS_{g,n}(W,G)(\gamma(1), \dotsc,\gamma(n)))$.
By linearity, it is extended to
\begin{equation*}
(\cvir)_{g,n} \colon \st^{\otimes n} \longrightarrow A^*(\sS_{g,n}(W,G)).
\end{equation*}

\subsection{PV virtual class on the $\Aut(W)$-invariant state space}
The evaluation of the virtual class on the states $e(\rR_{\gamma(1)}), \dotsc, e(\rR_{\gamma(n)})$ has a beautiful form.
In \cite[Section 2.4]{Guere1}, we find an explicit Koszul matrix factorization
\begin{equation*}
\kK(e(\rR_{\overline{\gamma}})) \in \textrm{MF}(\aA^{\overline{\gamma}},W_{\overline{\gamma}})
\end{equation*}
such that its Chern character is the element
\begin{equation*}
e(\rR_{\overline{\gamma}}) := e(\rR_{\gamma(1)}) \otimes \dotsc \otimes e(\rR_{\gamma(n)}).
\end{equation*}
Then, we reformulate Polishchuk and Vaintrob construction as follows.

We start with line bundles
\begin{equation}\label{fibreLR}
\cL^\rR_j:=\cL_j \Biggl(- \sum_{(\sigma_i,x_j) \in \rR_{\overline{\gamma}}} \sigma_i \Biggr)
\end{equation}
instead of $\cL_j$ and we apply the same procedure as Polishchuk and Vaintrob \cite[Sections 4.1-4.2]{Polish1}. to get resolutions by vector bundles
\begin{equation*}
R^\bullet\pi_* \cL^\rR_j = [A_j \rightarrow \widetilde{B}_j]
\end{equation*}
and morphisms\footnote{The morphism $\widetilde{\beta}_j$ comes from the resolution of $R^\bullet\pi_* \cL^\rR_j$ and the morphism $\widetilde{\alpha}_j$ arises from the algebraic relations \eqref{constraints}.}
\begin{equation}\label{widealpha}
\begin{array}{lcl}
\widetilde{\alpha}_j & \colon & \cO \rightarrow \cS^{a_{j+1}} A_{j+1}^\vee \otimes \widetilde{B}_j^\vee \oplus (\cS^{a_j-1} A_j^\vee \otimes A_{j-1}^\vee) \otimes \widetilde{B}_j^\vee, \\
\widetilde{\beta}_j & \colon & \widetilde{B}_j^\vee \rightarrow A_j^\vee.
\end{array}
\end{equation}

\noindent
Here, the convention is $(A_0,A_{N+1})=(0,A_N)$ for a chain polynomial and $(A_0,A_{N+1})=(A_N,A_1)$ for a loop polynomial.


At last, we get a two-periodic complex $(T,\delta)$ on the moduli space $\sS_{g,n}(W,G)(\overline{\gamma})$, given by infinite-rank vector bundles
\begin{equation*}
\begin{array}{lcl}
T^+ &:=& \mathrm{Sym} (A_1^\vee \oplus \dotsb \oplus A_N^\vee) \otimes \bigwedge_{\textrm{even}} (\widetilde{B}_1^\vee \oplus \dotsb \oplus \widetilde{B}_N^\vee), \\
T^- &:=& \mathrm{Sym} (A_1^\vee \oplus \dotsb \oplus A_N^\vee) \otimes \bigwedge_{\textrm{odd}} (\widetilde{B}_1^\vee \oplus \dotsb \oplus \widetilde{B}_N^\vee), \\
\end{array}
\end{equation*}
with the differential $\delta$ induced by \eqref{widealpha}.
By \cite[below Equation (34)]{Guere1} and \cite[Remark 1.5.1]{Polish1}, we have a quasi-isomorphism
\begin{equation*}
(T,\delta) \simeq p_*(\PV \otimes \kK(\rR_{\overline{\gamma}})).
\end{equation*}


\noindent
As a consequence, the virtual class evaluated at $e(\rR_{\overline{\gamma}})$ equals
\begin{equation}\label{virtualcomplex}
(\cvir)_{g,n}(e(\rR_{\overline{\gamma}})) = \Ch \left( H^+(T,\delta) - H^-(T,\delta) \right) ~ \prod_{j=1}^N \frac{\Td (\widetilde{B}_j)}{\Td (A_j)}.
\end{equation}

In genus zero and for chain polynomials\footnote{It works also with certain invertible polynomials, see \cite[Theorem 3.21]{Guere1} for a precise statement.}, the main result of \cite{Guere1} provides an explicit expression of the Chern character of the cohomology of $(T,\delta)$ in terms of the Chern characters of the higher push-forwards $R^\bullet\pi_*\cL^\rR_j$. The later are computed by Chiodo's formula \cite[]{Chiodo1}, using Grothendieck--Riemann--Roch theorem.
Thus the virtual class can be computed as well.

Interestingly, the same method provides an explicit computation of the cup product between the top Chern class of the Hodge bundle and the virtual class in arbitrary genus. We explain it in the following section.

\section{Polishchuk and Vaintrob's virtual class in higher genus}\label{contrib}
In this section, we prove our main theorem generalizing the computation of the virtual class in genus zero from \cite[Theorem 3.21]{Guere1} to Hodge integrals in arbitrary genus, see Theorem \ref{chainhigh}.

\subsection{Statement}
Let us consider an LG orbifold $(W,G)$ where $W$ is a Fermat monomial, a chain polynomial, or a loop polynomial and $G$ is an admissible group of symmetries.
We fix some elements $\gamma(1), \dotsc, \gamma(n) \in G$ and some admissible decorations $\rR_{\gamma(1)}, \dotsc, \rR_{\gamma(n)}$.
We consider the evaluation of the virtual class at the $\Aut(W)$-invariant state
\begin{equation*}
e(\rR_{\overline{\gamma}}) := e(\rR_{\gamma(1)}) \otimes \dotsc \otimes e(\rR_{\gamma(n)}).
\end{equation*}
In the case where $W$ is a loop polynomial, we further assume the existence of a variable $x_{j_0}$ such that
\begin{equation}\label{hyp}
\begin{array}{lcl}
\gamma_{j_0}(i) & \in & \langle e^{2\pi \ci \frac{w_{j_0}}{d}} \rangle ~~ \forall i, \\
w_{j_0} & \left| \right. & d, \\
\cL^\rR_{j_0} & = & \cL_{j_0} (-\sigma_1 - \dotsc - \sigma_n).
\end{array}
\end{equation}
By a cyclic permutation of the indices, we can assume that $j_0=N$.
Note also that conditions \eqref{hyp} are always true for a Fermat monomial or for the last variable $x_N$ of a chain polynomial.

\begin{thm}\label{chainhigh}
Let $(W,G)$ and $e(\rR_{\overline{\gamma}})$ be as above.
For any genus $g$, we have the following equality in the Chow ring of the moduli space of $(W,G)$-spin curves
\begin{equation}\label{formulelim3}
\begin{split}
\lambda_g^\vee ~ \cvirPV(e(\rR_{\overline{\gamma}}))_{g,n} & =  \lim_{t \to 1}  \prod_{j=1}^N \fc_{t_j}(-R^\bullet \pi_*(\cL^\rR_j)) \cdot \fc_{t_{N+1}}(\mathbb{E}^\vee) \\
 & =  \lim_{t \to 1}  \prod_{j=1}^N (1-t_j)^{r_j} \fc_{t_j}(-R^\bullet \pi_*(\cL_j)) \cdot \fc_{t_{N+1}}(\mathbb{E}^\vee), \\
\end{split}
\end{equation}
where $\lambda_g^\vee:=c_g(\mathbb{E}^\vee)$ is the top Chern class of the dual of the Hodge bundle, the integer $r_j := \mathrm{card} \left\lbrace i \left| \right. \gamma_j(i)=1\right\rbrace$ counts broad states, and
\begin{equation*}
t_{j+1} = \left\lbrace \begin{array}{ll}
t & \textrm{if $j=0$,} \\
t_j^{-a_j} & \textrm{if $1 \leq j \leq N-1$,} \\
t_N^{-d/w_N} & \textrm{if $j=N$.}
\end{array} \right. 
\end{equation*}
The characteristic class $\fc_t \colon K^0(S) \rightarrow A^*(S) [\![t]\!]$ is defined by
\begin{equation}\label{newclass}
\fc_t(B-A) = (1-t)^{-\Ch_0(A-B)} \exp \Biggl(\sum_{l \geq 1} s_l(t) \Ch_l(A-B) \Biggr),
\end{equation}
where the functions $s_l(t)$ are defined in \cite[Equation (67)]{Guere1} by
\begin{equation}\label{parametresl}
s_l(t) = \left\lbrace 
\begin{split}
&- \ln (1-t)  & \qquad \textrm{if $l=0$,} \\
&\cfrac{B_l(0)}{l} + (-1)^l \sum\limits_{k=1}^l (k-1)! \left( \frac{t}{1-t} \right)^k \gamma(l,k) & \qquad \textrm{if } l \geq 1, \\
\end{split}\right. 
\end{equation}
with the number $\gamma(l,k)$ defined by the generating function
\begin{equation*}
\sum_{l \geq 0} \gamma(l,k) \frac{z^l}{l!} := \frac{(e^z-1)^k}{k!}.
\end{equation*}
\end{thm}

Above theorem relies on our method developed in \cite[Section 3]{Guere1} together with two important observations:
\begin{itemize}
\item conditions \eqref{hyp} imply the algebraic relation
\begin{equation}
(\cL^\rR_N)^{\otimes \frac{d}{w_N}} \otimes \cO \hookrightarrow \omega_\cC,
\end{equation}
which is similar to relations \eqref{constraints},
\item the sheaf $\pi_* \omega$ is a vector bundle of rank $g$. It is called the Hodge bundle and we denote it by $\mathbb{E}$.
\end{itemize}
We now proceed to the proof of Theorem \ref{chainhigh}.

\subsection{Modified two-periodic complex and recursive complex}\label{appli}
The two above observations suggest us to introduce the line bundle
\begin{equation*}
\cL_{N+1} := \cO
\end{equation*}
and to choose a resolution $R^\bullet \pi_* \cL_{N+1} = [\cO \xrightarrow{0} \mathbb{E}^\vee]$ together with a morphism
\begin{equation}\label{newmorph}
\cO \rightarrow \cS^{d/w_N} A_N^\vee \otimes \mathbb{E}.
\end{equation}


Now, we consider the two-periodic complex $(\mathbf{T},\widetilde{\delta})$ with
\begin{eqnarray*}
\mathbf{T}^+ & = & \mathrm{Sym} (A_1^\vee \oplus \dotsb \oplus A_N^\vee) \otimes \Lambda_{\mathrm{even}} (\widetilde{B}_1^\vee \oplus \dotsb \oplus \widetilde{B}_N^\vee \oplus \mathbb{E}) \\
& = & T^+ \otimes \Lambda_{\mathrm{even}} \mathbb{E} \oplus T^- \otimes \Lambda_{\mathrm{odd}} \mathbb{E}
\end{eqnarray*}
and similarly for $\mathbf{T}^-$ exchanging odd and even, and with the differential
\begin{equation*}
\widetilde{\delta} = \delta_0 + \delta_1 + \delta_2,
\end{equation*}
where
\begin{itemize}
\item $\delta_0$ is induced by $\widetilde{\alpha}_1+\dotsc+\widetilde{\alpha}_{N-1}+\widetilde{\beta}_1+\dotsc+\widetilde{\beta}_N$,
\item $\delta_1$ is induced by $\widetilde{\alpha}_N$,
\item $\delta_2$ is induced by \eqref{newmorph}.
\end{itemize}
Note that the differential of the two-periodic complex $(T,\delta)$ is closely related to the differential $\delta_0+\delta_1$.

By the anticommutation relations among the maps $\widetilde{\alpha}_j, \widetilde{\beta}_j$ and \eqref{newmorph}, we obtain two double complexes
\begin{equation*}
(K_1=\mathbf{T},\delta_0 + \delta_1, \delta_2) \quad \textrm{and} \quad (K_2=\mathbf{T},\delta_0 + \delta_2, \delta_1).
\end{equation*}
The double complex $K_1$ is very explicit and we can write in particular
\begin{equation*}
(K_1)^{\pm,q} = T^\pm \otimes \Lambda^q \mathbb{E},
\end{equation*}
whereas the double complex $K_2$ is more involved.
Nevertheless, the cohomology groups of their associated two-periodic complexes agree and equal
\begin{equation*}
H^\pm(\mathbf{T},\delta_0 + \delta_1 + \delta_2).
\end{equation*}
We can abut to the total cohomology by looking at the spectral sequences given by the filtration by rows of these two double complexes.
In fact, the first page of the spectral sequence is even enough to compute the total cohomology in K-theory, as we show below.

On one side, we have
\begin{equation*}
(H^\pm(K_1,\delta_0 + \delta_1),\delta_2)^\bullet = (H^\pm(T,\delta) \otimes \Lambda^\bullet \mathbb{E},\delta_2),
\end{equation*}
which is a bounded complex of vector bundles by \cite[Equation (61)]{Guere1} and \cite[Theorem 3.3.1]{ChiodoJAG}, or \cite[Equation (1.20)]{Polish1}.
As a consequence, we have the following equalities in K-theory
\begin{eqnarray*}
H^+(\mathbf{T},\delta_0 + \delta_1 + \delta_2) & = & \bigoplus_{q \geq 0} (H^+(K_1,\delta_0 + \delta_1),\delta_2)^{2q} \oplus (H^-(K_1,\delta_0 + \delta_1),\delta_2)^{2q+1} \\
& = & H^+(T,\delta) \otimes \Lambda_{\mathrm{even}} \mathbb{E} \oplus H^-(T,\delta) \otimes \Lambda_{\mathrm{odd}} \mathbb{E}, \\
H^-(\mathbf{T},\delta_0 + \delta_1 + \delta_2) & = & H^+(T,\delta) \otimes \Lambda_{\mathrm{odd}} \mathbb{E} \oplus H^-(T,\delta) \otimes \Lambda_{\mathrm{even}} \mathbb{E}.
\end{eqnarray*}
Therefore, by the definition of the virtual class and by the equality
\begin{equation*}
\sum_{q \geq 0} (-1)^q \Ch(\Lambda^q V^\vee) \Td(V) = c_\textrm{top}(V)
\end{equation*}
for any vector bundle $V$, we obtain
\begin{equation}\label{membgauch}
\Ch (H^+(\mathbf{T},\widetilde{\delta})-H^-(\mathbf{T},\widetilde{\delta})) ~ \prod_{j=1}^N \frac{\Td(\widetilde{B}_j)}{\Td (A_j)} ~ \Td(\mathbb{E}^\vee) = \cvirPV(e(\rR_{\overline{\gamma}}))_{g,n} ~ c_\textrm{top}(\mathbb{E}^\vee).
\end{equation}

On the other side, we look at the cohomology groups
\begin{equation*}
H^\pm(K_2,\delta_0 + \delta_2).
\end{equation*}
The main point is that the two-periodic complex associated to $(K_2,\delta_0 + \delta_2)$ is a non-degenerate recursive complex with the vanishing condition\footnote{The vanishing condition comes from the fact that we can choose the resolution of $\mathbb{E}$ by vector bundles to be $\left[ 0 \rightarrow \mathbb{E}\right]$ since $\mathbb{E}$ is already a vector bundle.}, see \cite[Definitions 3.1,3.4 and Equation (40)]{Guere1}.
As a consequence, Theorem \cite[Theorem 3.5]{Guere1} implies that the cohomology groups are finite-rank vector bundles, so that
\begin{equation*}
H^+(K_2,\delta_0 + \delta_2)-H^-(K_2,\delta_0 + \delta_2) = H^+(\mathbf{T},\widetilde{\delta})-H^-(\mathbf{T},\widetilde{\delta}).
\end{equation*}
Furthermore, \cite[Theorem 3.19]{Guere1} yields an explicit computation of this difference in K-theory yielding
\begin{equation}\label{membdroit}
\Ch (H^+(\mathbf{T},\widetilde{\delta})-H^-(\mathbf{T},\widetilde{\delta})) ~ \prod_{j=1}^N \frac{\Td(\widetilde{B}_j)}{\Td (A_j)} ~ \Td(\mathbb{E}^\vee) = \lim_{t \to 1}  \prod_{j=1}^N \fc_{t_j}(-R^\bullet \pi_*(\cL^\rR_j)) \cdot \fc_{t_{N+1}}(\mathbb{E}^\vee),
\end{equation}
with $t_j$ and $\fc_t$ as in the statement of Theorem \ref{chainhigh}.
Equality between equations \eqref{membdroit} and \eqref{membgauch} proves the theorem.
\qed


\subsection{Some remarks}
Theorem \ref{chainhigh}, together with Chiodo's expression \cite[Theorem 1.1.1]{Chiodo1} of the Chern characters of $R^\bullet\pi_*\cL_j$ 
and Mumford's formula \cite[Equation (5.2)]{Mumford}, leads to explicit numerical computations of Hodge integrals that we have encoded into a MAPLE program \cites{computerprogram,PhDJG}.
Moreover, since the rank of the Hodge bundle is zero in genus zero, we easily recover \cite[Theorem 3.21]{Guere1}.



In particular, formula \eqref{formulelim3} gives some information in every genus on Polishchuk--Vaintrob virtual class for every Landau--Ginzburg orbifold $(W,G)$ with $W$ of chain type and $G$ any admissible group, provided that we evaluate the virtual class at $\Aut(W)$-invariant states.
In general, there are more broad states and we still need further work to understand how to deal with them (just as in genus zero).

In the generically semi-simple case, e.g.~when $G=\Aut(W)$, it is possible to compute the push-forward \eqref{pushfor} of the virtual class to the moduli space of stable curves in cohomology and in any genus using Teleman's result \cite{Teleman}.
Nevertheless, the answer is only in cohomology and is not on the virtual class itself but only on its push-forward \eqref{pushfor}.
Furthermore, it happens in general that FJRW theory is not generically semi-simple, e.g.~for $G=\mu_d$, and then Theorem \ref{chainhigh} is the very first systematic result for higher-genus virtual classes in this context.

An important application to the computation of Hodge integrals comes from \cite{Bur14}.
Indeed, Hodge integrals naturally appear in the definition of the double ramification hierarchy introduced by Buryak \cite{Bur14} and Theorem \ref{chainhigh} is then a useful tool to compute the equations of this integrable hierarchy.
Precisely, we wrote a specific computer program for $r$-spin theories \cites{computerprogram,PhDJG} and we proved a conjecture of Buryak when $r \leq 5$, see \cite[Theorem 1.1]{Guere3}.

At last, as already mentioned in the introduction, Theorem \ref{chainhigh} yields some tautological relations in the Chow ring of the moduli space of $(W,G)$-spin curves and therefore of the moduli space of stable curves.
Indeed, the right hand side of formula \eqref{formulelim3} is the limit of a power series with coefficients in the Chow ring of the moduli space of the theory.
We can develop it and express it as a Laurent series in $\epsilon := t^{-1}-1$ to find an expression like
\begin{equation*}
C_{-p} \cdot \frac{1}{\epsilon^p} + C_{-p+1} \cdot \frac{1}{\epsilon^{p-1}} + \dotsb + C_{-1} \cdot \frac{1}{\epsilon} + C_0 + C_1 \cdot \epsilon + \dotsb.
\end{equation*}
According to the discussion before \cite[Corollary 3.20]{Guere1}, above expression has the property that
\begin{equation*}
C_m \in \bigoplus_{k \geq \mathrm{degvir}+g -m} A^{k}(\sS_{g,n}(W,G)(\overline{\gamma}))
\end{equation*}
and $p=2g-3+n-\mathrm{degvir}$, where the integer $\textrm{degvir}:=-\sum_j\Ch_0(R^\bullet\pi_* \cL^\rR_j)$ is the Chow degree of the virtual class.
As a consequence of the existence of the limit in \eqref{formulelim3} when $\epsilon \to 0$, we obtain relations
\begin{equation*}
C_m = 0 \textrm{ for } m<0.
\end{equation*}
In a work in progress \cite{team4}, we compare the push-forward to $\overline{\cM}_{g,n}$ of these relations with Pixton's relations, see \cite{Pixton}.

\bibliographystyle{plain} 
\bibliography{bibliothese} 


\end{document}